\input amstex
\documentstyle{amsppt}
\hsize=15.34 truecm
\vsize=22.83 truecm

\topmatter
\topinsert
\captionwidth{10 truecm}
\flushpar Chin.
Sci. Bulletin, 1999, 44(23), 2465--2470 (Chinese Ed.); 2000, 45:9
(English Ed.), 769--774
\endinsert
\rightheadtext{Mu-Fa Chen}
\leftheadtext{Eigenvalues, inequalities and ergodic theory}

\endtopmatter
\font\MR=cmbx10 scaled \magstephalf

\font\cmt=cmss9
\newsymbol\leqs 1336
\newsymbol\geqs 133E

\document

\def\d{\text{\rm d}}

\def\var{\text{\rm Var}}

\centerline{\MR Eigenvalues, inequalities and ergodic theory}

\smallskip

\centerline{Mu-Fa Chen}

\centerline{(Beijing Normal University, Beijing 100875)}

\medskip

\flushpar{\bf Abstract}\quad{\cmt This paper surveys the main
results obtained during the period 1992--1999 on three aspects
mentioned at the title. The first result is a new and general
variational formula for the lower bound of spectral gap (i.e., the
first non-trivial eigenvalue) of elliptic operators in Euclidean
space, Laplacian on Riemannian manifolds or Markov chains (\S 1).
Here, a probabilistic method---coupling method is adopted. The new
formula is a dual of the classical variational formula. The last
formula is actually equivalent to Poincar\'e inequality. To which,
there are closely related logarithmic Sobolev inequality, Nash
inequality, Liggett inequality and so on. These inequalities are
treated in a unified way by using Cheeger's method which comes
from Riemannian geometry. This consists of \S 2. The results on
these two aspects are mainly completed by the author joint with F.
Y. Wang. Furthermore, a diagram of the inequalities and the
traditional three types of ergodicity is presented (\S 3). The
diagram extends the ergodic theory of Markov processes. The
details of the methods used in the paper will be explained in a
subsequent paper under the same title.}

\flushpar {\bf  Keywords}\quad Eigenvalue\quad inequality\quad
ergodic theory \quad Markov process

\medskip

\flushpar{\MR 1\quad New variational formula for the lower bound
of spectral gap}

\smallskip

\flushpar{\bf 1.1 \quad Story of estimating $\lambda_1$ in
geometry}

We recall the study on
$ \lambda_1$ in geometry. From the story below, one should have
some feeling about the difficulty of the hard mathematical topic.

Consider Laplacian $ \Delta $ on a compact Riemannian manifold
$(M, g)$, where $g$ is the Riemannian metric. The spectrum of $
\Delta $ is discrete: $\cdots \leqs - \lambda_2\leqs -\lambda_1 <
- \lambda_0=0$ (may be repeated). Estimating these eigenvalues $
\lambda_k$ (especially $ \lambda_1$) consists an important section
and chapter of the modern geometry. As far as we know, until now,
five books have been devoted to this topic. Here we list only the
geometric books but ignore the ones on general spectral
theory$^{\text{[1]---[5]}}$. Denote by $d$, $D$ and $K$
respectively the dimension, the diameter and the lower bound of
Ricci curvature (Ricci$_M \geqs K g$) of the manifold $M$. We are
interested in estimating $ \lambda_1$ in terms of these three
geometric quantities. For an upper bound, it is relatively easy.
Applying a test function $f\in C^1(M)$ to the classical
variational formula
$$\lambda_1=\inf\{\tsize
\int_M \| \nabla f\|^2 : \;  f \in  C^1 (M), \; \tsize\int f\d x=0,
\; \tsize\int f^2\d x =1 \},$$
where ``$\d x$'' is the Riemannian volume element, one gets an
upper bound. However,  the lower bound is much harder.  The
previous works have studied the lower estimates case by case  by
using different elegant methods. Eight of the most beautiful lower
bounds are listed in the following table.
$$\matrix
\format\l& \l& \l& \l\\
&\text{A. Lichnerowicz (1958)} &\qquad\dfrac{d}{d-1}\, K,
\quad K\geqs 0. &\qquad (1)\\
&{}&{}&{}\\
&{\matrix\text{P. H. B\'erard, G. Besson}\\
\text{ \& S. Gallot (1985)}\endmatrix} &\qquad d\,  \bigg\{
\dfrac{\int_0^{\pi/2}\cos^{d-1}t \d t}{
\int_0^{D/2}\cos^{d-1}t \d t }\bigg\}^{2/d}, \quad K=d-1>0. &\qquad (2)\\
&{}&{}&{}\\
&{\matrix\text{P. Li \& S. T. Yau (1980)}\\
\text{   }\endmatrix}&\qquad \dfrac{\pi^2}{2\,D^2},\quad K\geqs 0.
&\qquad  (3)\\
&{}&{}&{}\\
&{\matrix\text{J. Q. Zhong \& H. C. Yang (1984)}\\
\text{   }\endmatrix}&\qquad\dfrac{\pi^2}{D^2},
\quad K\geqs 0. &\qquad (4)\\
&{}&{}&{}\\
&{\matrix\text{P. Li \& S. T. Yau (1980)}\\
\text{  }\endmatrix}&\qquad\dfrac{1}{D^2 (d - 1)\exp\big[1 +
 \sqrt{1 + 16\alpha^2}\big]},    \quad K\leqs 0. &\qquad (5) \\
&{}&{}&{}\\
&{\matrix\text{K. R. Cai (1991)}\\
\text{   }\endmatrix}&\qquad\dfrac{\pi^2}{D^2} + K,
\quad K\leqs 0.  &\qquad  (6)\\
&{}&{}&{}\\
&{\matrix\text{H. C. Yang (1989) \& F. Jia (1991)}\\
\text{  }\endmatrix}&\qquad \dfrac{\pi^2}{D^2}
e^{-\alpha},\quad\text{if }\;d\geqs 5,
          \quad K\leqs 0.  &\qquad (7)\\
&{}&{}&{}\\
&{\matrix\text{H. C. Yang (1989) \& F. Jia (1991)}\\
\text{  }\endmatrix}&\qquad \dfrac{\pi^2}{2\,D^2} e^{-\alpha'},
\quad \text{if }\; 2\leqs d\leqs 4, \quad K\leqs 0, &\qquad (8)
\endmatrix$$

\flushpar where $\alpha=D\sqrt{|K|(d-1)}/2,\;\; \alpha' =
D\sqrt{|K|((d-1)\vee 2)}/2$. All together, there are five sharp
estimates ((1), (2), (4), (6) and (7)). The first two are sharp
for the unit sphere in two- or higher-dimension but it fails for
the unit circle; the fourth, the sixth and the seventh estimates
are all sharp for the unit circle. The above authors include
several famous geometers and the estimates were awarded several
times. From the table, it follows that the picture is now very
complete, due to the effort by the geometers in the past 40 years.
For such a well-developed field, what can we do now? Our original
starting point is to learn from the geometers, study their
methods, especially the recent new developments. It is surprising
that we actually went to the opposite direction, that is, studying
the first eigenvalue by using a probabilistic method. It was
indeed not dreamed that we could finally find a general formula.

\flushpar {\bf 1.2\quad New variational formula}

To state the result, we need two notations
$$\align
&C(r)=\text{cosh}^{d-1}\bigg[ \dfrac{r}{2}\sqrt{ \dfrac{-K}{d-1}}\bigg],
           \quad\quad r\in (0, D).\\
&{\Cal F}=\{f\in C[0,D]: f>0 \text{ on } (0,D)\}.\endalign$$ Here
the dimension $d$,  the diameter $D$ and the lower bound of Ricci
curvature $K$ have all been used.

{\bf Theorem [General formula]\,}(Chen \& Wang$^{[6]}$).
$\lambda_1\geqs \sup\limits_{f\in {\Cal F}}\inf\limits_{r\in
(0,D)} \dfrac{4 f(r)}{\int_0^r  C(s)^{-1}\d s\int_s^D C(u)f(u)\d
u}$.

The new variational formula has its essential value in estimating
the lower bound. It is a dual of the classical variational formula
in the sense that ``$\inf$'' is replaced by ``$\sup$''. The last
formula  goes back to Lord S. J. W. Rayleigh(1877) or
E.\,Fischer\,(1905). Noticing that there are no common points in
these two formulas, this explains the reason why such a formula
never appeared before. Certainly, the new formula can produce a
lot of new lower bounds. For instance, the one corresponding to
the trivial function $f{\Cal E}uiv 1$ is still non-trivial in
geometry. Next, let $\alpha$ be the same as above and let
$\beta=\dfrac{\pi}{2D}$. Applying the formula to the test
functions $\sin(\beta r)$, $\sin(\alpha r)$, $\sin(\beta r)$ and
$\cosh^{d-1}(\alpha r)\sin(\beta r)$ successively, we obtain the
following:

 {\bf Corollary \,}(Chen\,\&\,Wang$^{[6]}$).
$$\align
&\lambda_1\geqs\frac{\pi^2}{D^2}+\max\Big\{\frac{\pi}{4d},1-\frac{2}{\pi}\Big\}K,
  \qquad K\geqs 0 \tag 9\\
& \lambda_1 \geqs \frac{dK}{d-1}\bigg\{1-\cos^d\bigg[\frac{D}{2}
 \sqrt{\frac{K}{d-1}}\bigg]\bigg\}^{-1} , \qquad d>1, \quad K\geqs 0 \tag 10\\
& \lambda_1\geqs \frac{\pi^2}{D^2} +\Big(\frac{\pi}{2}-1\Big)K,\qquad K\leqs 0 \tag 11\\
&\lambda_1\geqs
\frac{\pi^2}{D^2}\sqrt{1-\frac{2D^2K}{\pi^4}}\cosh^{1-d}
\bigg[\frac{D}{2} \sqrt{\frac{-K}{d-1}}\bigg] \qquad \qquad d>1,
\quad K\leqs 0. \tag 12
\endalign$$

 {\bf Comments}. \roster
\item The corollary improves all the estimates (1)---(8).
(9) improves (4); (10) improves (1) and (2); (11) improves
(6); (12) improves (7) and (8).
\item The theorem and corollary valid also for the manifolds with
convex boundary with Neumann boundary condition. In this case,
the estimates (1)---(8) are believed by geometers to be true.
However, only the Lichnerowicz's estimate (1) was
proved by J. F. Escobar until 1990. Except this, the others in (2)---(8)
(and furthermore  (9)---(12)) are all new in geometry$^{[6]}$.
\item For more general non-compact manifolds, elliptic operators or
Markov chains, we also have the corresponding dual variational
formula$^{[7], [8]}$. The point is that only three parameters $d$ , $D$ and $K$
are used in the geometric case, but there are infinite parameters
in the case of elliptic operators or Markov chains. Thus, the latter
cases are more complicated. Actually, the above formula is a
particular example of our general formula for elliptic operators.
In dimensional one, our formula is complete.
\item The probabilistic method---coupling method was developed
by the present author before this work for more than ten years.
The above study  was the
first time for applying the method to estimating the eigenvalues. For
a long time, almost nobody believes that the method can achieve
sharp estimate. From these facts, the influence of the
above results to probability theory and spectral theory should
be clear$^{[8]}$.\endroster

\medskip

\flushpar{\MR 2\quad Basic inequalities and new forms of Cheeger's
constants}
\smallskip

\flushpar {\bf 2.1\quad Basic inequalities}

Let $(E, \Cal E , \pi )$ be a probability space satisfying $\{(x,
x): x\in E\}\in {\Cal E} \times {\Cal E}$. Denote by $L^p(\pi)$
the usual real $L^p$-space with norm $\| \cdot \|_p$. Write
$\|\cdot\|=\|\cdot\|_2$. Our main object is a symmetric form  $(D,
\Cal D (D))$ on $L^2(\pi)$. For Laplacian on manifold, the form
used in the last part is the following
$$D(f):=D(f, f)=\int_M \|\nabla f\|^2 \d x, \quad\quad \Cal D (D)\supset C^ \infty
(M).$$
Here, only the diagonal elements $D(f)$ is written, but the non-diagonal
elements can be then deduced from the diagonal ones by using the
quadrilateral role. The classical variational formula for spectral gap
now can be rewritten into the following form.
$$\text{{\it Poincar\'e inequality}}: \qquad
\var(f)\leqs C D(f), \quad\quad f\in L^2(\pi)$$
where $\var(f)=\pi(f^2)-\pi(f)^2$, $\pi(f)= \int f\d \pi$ and
$C(= \lambda_1^{-1})$
is a constant. Thus, the study on the spectral gap is the same as the one on
Poincar\'e inequality of the form $(D, \Cal D (D))$. Nevertheless, we have more
symmetric forms. For an elliptic operator in $\Bbb R^d$, the corresponding form
is as follows.
$$D(f)=\frac{1}{2}\int_{\Bbb R^d}\langle a(x)\nabla f(x), \nabla
f(x)\rangle\pi (\d x), \qquad \Cal D(D)\supset C_0^ \infty (\Bbb
R^d),$$ where $\langle \cdot, \cdot \rangle$ denotes the standard
inner product in $\Bbb R^d$ and $a(x)$ is positive definite.
Corresponding to an integral operator (or symmetric kernel) on
$(E, \Cal E )$, we have the symmetric form
$$D(f)= \frac{1}{2}\int_{E\times E}J(\d x, \d y)[f(y)-f(x)]^2,\qquad
\Cal D (D)=\{f\in L^2(\pi): D(f)< \infty\}, \tag 13$$
where $J$ is a non-negative, symmetric measure having no charge on
the diagonal set $\{(x, x): x\in E\}$. A typical example in our
mind is the reversible jump process with $q$-pair $(q(x),$ $ q(x, \d y))$
and reversible measure $\pi$.
Then $J(\d x, \d y)= \pi(\d x)q(x, \d y)$. More especially, for a
reversible $Q$-matrix $Q=(q_{ij})$ with reversible measure
$(\pi_i>0)$, we have density $J_{ij}=\pi_i q_{ij}\,(j\ne i)$ with
respect to the counting measure.

For a given symmetric form $(D, \Cal D(D) )$, except Poincar\'e inequality,
there are also other basic inequalities.
$$\align
&\text{{\it Nash inequality}}:\;\;\;
  \quad\quad  \var(f)\leqs C D(f)^{1/p} \|f\|_1^{2/q}, \quad\quad f\in L^2(\pi)\\
&\text{{\it Liggett inequality}}: \quad\quad  \var(f)\leqs C
D(f)^{1/p} \text{Lip}(f)^{2/q}, \quad\quad f\in L^2(\pi)
\endalign$$ where $C$ is a constant and Lip$(f)$ is the Lipschitz
constant of $f$ with respect to some distance $\rho$. The above
three inequalities are actually particular cases of the following
one
$$\text{{\it Liggett-Stroock inequality}}: \quad\quad
\var(f)\leqs C D(f)^{1/p} V(f)^{1/q}, \quad\quad f\in L^2(\pi)$$
where $V: L^2(\pi)\to [0, \infty]$ is homogeneous of degree two:
$V(c_1 f+c_2)=c_1^2 V(F)$,
$c_1$, $c_2\in \Bbb R$. Another closely related one is
$$\text{{\it Logarithmic Sobolev inequality}}: \qquad
\int f^2\log \big( f^2/\ |f\|^2\big) \d\pi \leqs  C D(f), \quad\quad f\in L^2(\pi).$$
\smallskip

\flushpar {\bf 2.2\quad Statue of the research}

From now on, we
restrict ourselves to the symmetric form (13) corresponding to integral
operators. The question is under what condition on the symmetric
measure $J$, the above inequalities hold.
In contrast with the probabilistic method used in the last part,
here we adopt Cheeger's method (1970) which comes from
Riemannian geometry.

We call $\lambda _1:=\inf\{D(f): \pi (f)=0,\, \|f\|=1\}$ the
{\it spectral gap} of the form $(D, \Cal D (D))$. For bounded
jump processes, the main known result is the following.
\smallskip

{\bf Theorem} (Lawler \& Sokal (1988)). $\lambda_1\geqs
\dfrac{k^2}{2M}$, { where} $k=\inf\limits_{\pi (A)\in (0, 1)}
\dfrac{\int_A \pi(\d x) q(x, A^c)}{\pi(A)\wedge \pi(A^c)}$,
$M=\sup\limits_{x\in E} q(x)$.
\smallskip

In the past seven years, the theorem has been collected into
six books$^{\text{[9]---[14]}}$.
From the titles of the books, one sees the wider range of the
applications of the study. The problem is:
the result fails for unbounded operator. Thus, it has been a challenge
open problem in the past ten years or more to handle the unbounded
situation.

As for logarithmic Sobolev inequality, there is a large number of
publications in the past twenty years or more for differential operators.
However, there was almost no result for integral operators
until the next result appeared.
\smallskip

{\bf  Theorem} (Diaconis \& Saloff-Coste (1996)). {  Let $E$
be a finite set and $\sum_{j} |q_{ij}|=1$ holds for all $i$. Then
the {\it logarithmic Sobolev constant} $ \sigma:= \inf
\big\{D(f)/\int f^2 \log[|f|/\|f\|]: \|f\|=1 \big \}$ satisfies
$\sigma\geqs \dfrac{2(1-2\pi_*)\lambda_1}{\log[1/\pi_*-1]}$, where
$\pi_*=\min_i \pi_i$.}
\smallskip

Obviously, the result fails again for infinite $E$. The problem is due to
the limitation of the method used in the proof.
\smallskip

\flushpar {\bf 2.3\quad New result}

Corresponding to  three inequalities,
we introduce respectively the following new forms of Cheeger's constants.
$$\matrix
\format\l\qquad &\l\\
\underline{\text{Inequality}}&\qquad\qquad\qquad\underline{\text{Constant }k^{(\alpha)}}\\
{ }& { }\\
\text{Poincar\'e}  &\inf\limits_{\pi(A)\in (0, 1)}
\dfrac{J^{(\alpha)}(A\times A^c)}{\pi (A)\wedge \pi(A^c)}\qquad
{(\text{Chen \& Wang}^{[15]})}\\
{ }&{ }\\
\text{Nash} &
 \inf\limits_{\pi(A)\in (0, 1)}\dfrac{J^{(\alpha)}(A\times A^c)}
 {[\pi (A)\wedge \pi(A^c)]^{(\nu-1)/\nu}}\qquad {(\text{Chen}^{[16]})}\\
{ } &  \nu=2(q-1)\\
\text{Log. Sobolev} &
\lim\limits_{r\to 0}\inf\limits_{\pi(A)\in(0, r]}
 \dfrac{J^{(\alpha)}(A\times A^c)}
{\pi (A) \sqrt{\log[e +  \pi(A)^{-1}]}}\qquad {(\text{Wang}^{[17]})}\\
\text{  } & \lim\limits_{\delta\to \infty} \inf\limits_{\pi(A)>0}
\dfrac{J^{(\alpha)}(A\times A^c) + \delta\pi(A)} {\pi (A) \sqrt{1-
\log \pi(A)}}\qquad {(\text{Chen}^{[18]})}\endmatrix$$ where $r(x,
y)$ is a symmetric, non-negative function such that $J^{( \alpha
)}(\d x, \d y)
   := I_{\{r(x, y)>0\}}\dfrac{J(\d x, \d y)}{r(x, y)^{ \alpha }}$
$(\alpha >0)$ satisfies $\dfrac{J^{(1)}(\d x, E)}{\pi(\d x)}\leqs
1$, $\pi$-a.s. For convenience, we use the convention $J^{(0)}=J$.
Now, our main result can be easily stated as follows.
\smallskip

{\bf  Theorem}. $k^{(1/2)}>0\Longrightarrow$
{ the corresponding inequality holds}.
\smallskip

The result is proved in four papers [15]---[18].
At the same time, some estimates for the upper or lower bounds
are also presented. These estimates can be sharp or qualitatively sharp,
which did not happen before in using Cheeger's technique.
\medskip

\flushpar{\MR 3 \quad New picture of ergodic theory}
\smallskip

\flushpar {\bf 3.1\quad Importance of the inequalities}

Let $(P_t)_{t\geqs 0}$
be the semigroup determined by the symmetric form
$(D, \Cal D (D))$. Then, various applications of the inequalities are
based on the following result.
\smallskip

{\bf  Theorem}. {
\roster
\item Let $V(P_t f)\leqs V(f)$ for all $t\geqs 0$ and $f\in L^2(\pi)$ (which is
automatic when $V(f)=\|f\|_r^2$). Then Liggett-Stroock inequality implies that
$$\var(P_t f) \leqs C V(f)/t^{q-1}, \qquad t> 0. \tag 14$$
\item Conversely, $(14)\Longrightarrow $ Liggett-Stroock inequality.
\item Poincar\'e inequality $\Longleftrightarrow \var(P_t f)\leqs \var(f)\exp[-2 \lambda_1 t]$.
\endroster}
\smallskip

Note that $\var (P_tf)=\|P_t f -\pi (f)\|^2$. Therefore, the above inequalities describe
some type of $L^2$-ergodicity of the semigroup $(P_t)_{t\geqs 0}$.
In particular, we call (14) $L^2$-{\it algebraic convergence}.
These inequalities have become powerful tools in
the study on infinite-dimensional mathematics (phase transitions,
for instance) and the effectiveness of random algorithms.
\smallskip

\flushpar {\bf 3.2\quad Three traditional types of ergodicity}

In the study of  Markov
processes, the following three types of ergodicity are well known.
$$\align
&\text{\it Ordinary ergodicity}:\qquad \lim_{t\to  \infty}\|p_t(x, \cdot)-\pi\|_{\var}=0\\
&\text{\it Exponential ergodicity}: \qquad \|p_t(x, \cdot)-\pi\|_{\var}\leqs C(x) e^{- \varepsilon t}\\
&\text{\it Strong ergodicity}:\qquad\;
  \lim_{t\to  \infty}\sup_x \|p_t(x, \cdot)-\pi\|_{\var}=0\endalign$$
where $p_t(x, \d y)$ is the transition function of the Markov process and
$\|\cdot\|_{\var}$ is the total variation norm. They obey the following relation:
$\text{Strong ergodicity}\Longrightarrow\text{Exponential ergodicity}
\Longrightarrow\text{Ordinary ergodicity}$. Now, it is natural to ask
the following question. Does there exist any relation between the above
inequalities and the traditional three types of ergodicity?
\smallskip

\flushpar{\bf 3.3\quad New picture of ergodic theory}
\smallskip

{\bf Theorem}$^{[16], [19], [20]}$. {For reversible Markov chains,
we have the following diagram:
$$\matrix
\format\c\text{\hskip-3em}& \c\text{\hskip-3em} &\c\\
{ }& \text{Nash inequality}&{ }\\
\text{\hskip4em}\swarrow \! \! \! \!\swarrow&{ }
 &{\text{\hskip-4em}\searrow \! \! \! \! \searrow }\\
\text{Log. Sobolev inequality}&{ }&\text{Strong ergodicity}\\
\Downarrow &{  } &\Downarrow\\
\text{Poincar\'e  inequality\quad}&\Longleftrightarrow&
    \text{\quad exponential ergodicity} \\
{ }&\Downarrow&{ }\\
{ }&\text{$L^2$-algebraic ergodicity}&{ }\\
{ }&\Downarrow&{ }\\
{ }& \text{Ordinary ergodicity}&{ }
\endmatrix$$
where $L^2$-{\it algebraic ergodicity} means that (14) holds for some
$V$ having the properties: $V$ is homogeneous of degree two, $V(f)<\infty$
for all functions $f$ with finite support.   }
\smallskip

{\bf  Comments}.\roster
\item The diagram is complete in the following sense. Each single-side
implication can not be replaced by double-sides one. Moreover, strong ergodicity
and logarithmic Sobolev inequality are not comparable.
\item The application of the diagram is obvious. For instance,
one obtains immediately some criteria (which are indeed new)
for Poincar\'e inequality to be held
from the well-known criteria for the exponential ergodicity. On the
other hand, by using the estimates obtained from the study on
Poincar\'e inequality, one may estimate exponentially ergodic
convergence rate (for which, the knowledge is still very limited).
\item Except the equivalence, all the implications in the
diagram are suitable for more general Markov processes.
The equivalence in the diagram should be also suitable for more
Markov processes but it may be false in the infinite-dimensional situation.
\item No doubt, the diagram extends the ergodic theory of Markov processes.
\endroster
\smallskip

{\bf Acknowledgement}\quad {Research supported in part by NSFC
(No. 19631060), Math. Tian Yuan Found., Qiu Shi Sci. \& Tech.
Found., RFDP and MCME.}

\medskip

\centerline{\bf References}
\widestnumber\no {10000}
\ref\no [1]
\by Chavel I.
\book {\it Eigenvalues in Riemannian Geometry}
\publ New York: Academic Press, 1984 \endref

\ref\no [2]
\by B\'erard P H
\book {\it Spectral Geometry: Direct and Inverse Problem}
\publ {\bf LNM.} vol 1207, New York: Springer-Verlag, 1986\endref

\ref\no [3]
\by Yau S T, Schoen R
\book Differential Geometry $($In Chinese$)$
\publ Beijing: Science Press,   1988\endref

\ref\no [4]
\by Li P
\book {\it Lecture Notes on Geometric Analysis}
\publ Seoul National Univ,  Korea, 1993\endref

\ref\no [5]
\by Ma C Y
\book The Spectrum of Riemannian Manifolds $($In Chinese$)$
\publ  Nanjing: Press of Nanjing U,  1993\endref

\ref\no [6]
\by Chen M F, Wang F Y
\paper General formula for lower bound of the first
      eigenvalue
\jour Sci Sin, 1997, 40:4, 384--394
\endref

\ref\no [7]
\by Chen M F, Wang F Y
\paper  Estimation of spectral gap for elliptic operators
\jour Trans Amer Math Soc, 1997, 349: 1239--1267\endref

\ref\key [8]
\by Chen M F
\paper Coupling, spectral gap and related topics
\jour Chin Sci Bulletin, 1997, (I):
     42:16, 1321--1327; (I\!I):
     42:17, 1409--1416;
(I\!I\!I):
 42:18, 1497--1505
\endref

\ref\no [9]
\by Chen M F
\book From Markov  Chains to Non-Equilibrium Particle Systems
\publ Singapore: World Scientific, 1992 \endref

\ref\no [10]
\by Sinclair A
\book Algorithms for Random Generation and Counting: A  Markov
Chain Approach
\publ Boston: Birkh\"auser, 1993\endref

\ref\no [11]
\by Saloff-Coste L
\paper Lectures on finite Markov chains
\jour LNM {\bf 1665}, 301--413, New York: Sprin\-ger-Verlag,
1997\endref

\ref\no [12]
\by Chung F R K
\book Spectral Graph Theory
\publ CBMS, {\bf 92}, Rhode Island: AMS, Providence,  1997
\endref

\ref\no [13]
\by Colin de Verdi\`ere Y
\book Spectres de Graphes
\publ Paris: Publ Soc Math France, 1998\endref

\ref\no [14]
\by Aldous D G \& Fill J A (1994--)
\book Reversible Markov Chains and Random Walks on Graphs
\publ URL {\bf www. stat.Berkeley.edu/users/aldous/book.html}\endref

\ref\key [15] \by Chen M F, Wang F Y \paper Cheeger's inequalities
for general symmetric forms and existence criteria for spectral
gap \jour {\bf Abstract}. Chin Sci Bulletin, 1998, 43:18,
1516--1519. Ann. Prob. 2000, 28:1, 235--257
\endref

\ref\key [16]
\by Chen M F
\paper Nash inequalities for general symmetric forms
\jour Acta Math Sin Eng Ser, 1999, 15:3, 353--370\endref

\ref\key [17]
\by Wang F Y
\paper Sobolev type inequalities for general symmetric forms
\jour to appear in Proc Amer Math Soc, 1999\endref

\ref\key [18] \by Chen M F \paper Logarithmic Sobolev inequality
for symmetric forms \jour Sci Chin, 2000, 43:6, 601--608\endref

\ref\key [19]
\by Chen M F
\paper Equivalence of exponential ergodicity and $L^2$-exponential
  convergence for Mar\-kov chains
\jour Stoch Proc Appl, 2000,87, 281--297 \endref

\ref\key [20]
\by Chen M F
\paper A new story of ergodic theory
\jour to appear in Proceedings of IMS Workshop on Applied Probability,
Hong Kong: Intern. Press \endref

\smallskip

{\hskip30em (Received August 2, 1999)}
\enddocument